\documentclass[11pt,reqno]{amsart}
\usepackage{amscd,amssymb,amsmath,amsfonts,latexsym,amsthm}
\usepackage[T1]{fontenc}
\usepackage{inputenc}
\usepackage{mathtools}
\usepackage[arrow,matrix]{xy}
\usepackage{graphicx}
\usepackage{cite}
\usepackage{ifpdf}
\usepackage{tikz}
\usepackage{longtable}
\usepackage{setspace}

\usepackage{hyperref}

\oddsidemargin=0.3in \evensidemargin=0.3in \theoremstyle{plain}
\newtheorem{thm}[subsection]{Theorem}

\newtheorem{rk}[subsection]{Remark}
\newtheorem{defn}[subsection]{Definition}

\numberwithin{equation}{section} 

\newcommand{\F}{{\mathcal F}}
\newcommand{\M}{{\mathcal M}}

\newcommand{\bea}{\begin{eqnarray}}
\newcommand{\eea}{\end{eqnarray}}

\newcommand{\C}{\mathbb{C}}

\begin{document}
\title[On chains and Rota-Baxter operators of evolution algebras.]
{On chains and Rota-Baxter operators of evolution algebras}

\author{M. Ladra, Sh.N. Murodov}

\address{M. \ Ladra\\ Department of Mathematics \& Institute of Mathematics, University of Santiago de Compostela, Santiago de Compostela, Spain.} \email {manuel.ladra@usc.es}
\address{Sh.\ N.\ Murodov\\ Bukhara State Medical institute,
	Bukhara, Uzbekistan \& Institute of Mathematics, University of Santiago de Compostela, Santiago de Compostela, Spain.} \email {murodovs@yandex.ru}

\begin{abstract} The paper is devoted to study new classes of chains of evolution algebras and
	their time-depending dynamics. Moreover, we construct some Rote-Baxter operators of such algebras.
\end{abstract}


\maketitle

\section{Introduction} %\label{sec:intro}

In the 1920s and 1930s a new object - general genetic algebra was introduced to mathematics,
which was the product of interactions between Mendelian genetics
and mathematics.
Recognizing algebraic structures and properties in
Mendelian genetics was one of the most important steps to start study genetic algebras.
 Mendel first exploited some symbols \cite{men}, which is quite
algebraically suggestive to express his genetic laws. In fact, it was later termed
``Mendelian algebras'' by several authors. Mendel's laws were mathematically formulated by Serebrowsky \cite{ser},
who was the first to give an
algebraic interpretation of the sign ``$\times$'', which indicated sexual reproduction. Later Glivenkov \cite{gli}
introduced the so-called Mendelian
algebras. Independently, Kostitzin \cite{kos} also introduced a ``symbolic multiplication'' to express
Mendel's laws. Etherington \cite{e1,e2,e3} gave the systematic study of algebras occurring in genetics
and  introduced the formal language of abstract algebra to the
study of genetics. These algebras, in general, are non-associative.

There exist several classes of non-associative algebras (baric, evolution,
stochastic, etc.), whose investigation has provided a number of significant
contributions to theoretical population genetics. Such classes have been
defined different times by several authors, and all algebras belonging to
these classes are generally called ``genetic''.

Fundamental contributions have also been
made by Gonshor \cite{gon}, Schafer \cite{sch}, Holgate~\cite{hol1,hol2}, Heuch \cite{heuch}, Reiers\"ol \cite{reiers},
Abraham \cite{abr}. Until the 1980s, the most comprehensive reference in this
area was W\"orz-Busekros' book \cite{w}. More recent results, such as evolution theory in genetic algebras, can be found in Lyubich's book \cite{ly}. A good survey
article is Reed's paper \cite{reed}. All algebras studied by these authors are generally called ``genetic''.

 In the present day, Non-Mendelian genetics is a basic language of molecular geneticists.
 Naturally, the question arises: what new subject non-Mendelian genetics offers to mathematics?
What mathematics offers to understanding of non-Mendelian genetics?
The systematic formulation of reproduction in non-Mendelian genetics as
multiplication in algebras was introduced in \cite{jpp} and called as ``evolution algebras''.
These are algebras in which the multiplication tables are motivated by evolution laws of genetics.

 The foundation of evolution algebra theory and
applications in non-Mendelian genetics and Markov chains are
developed by Tian in \cite{t}. Evolution algebras associated to function spaces defined by graphs, state spaces, and Gibbs measures are studied in \cite{rt}.

A notion of chain of evolution algebras was introduced in \cite{clr}, where the sequence
of matrices of  structural constants of the chain of evolution
algebras satisfies an analogue of the Chapman-Kolmogorov equation.

Originally Rota-Baxter operators were defined on associative algebras by G. Baxter to
solve a problem in probability \cite{bax} and then developed by the Rota school \cite{rota}.
These operators have connections with many areas of mathematics and mathematical physics,
such as number theory, combinatorics, quantum field theory.

The paper is organized as follows. In Section~\ref{S:chain}  we give the main definitions related to a chain of evolution algebras.
In Section~\ref{S:rota} we construct new chains of evolution algebras (CEAs), and study theirs time-depending dynamics.

Finally, in Section~\ref{S:rota} we give Rota-Baxter operators on two-dimensional complex evolution algebras.

\section{Chain of Evolution Algebras} \label{S:chain}

Evolution algebras are defined as follows.

\begin{defn}
	Let $(E,\cdot)$ be an algebra over a field $K$. If it admits a
	basis $\{e_1,e_2,\dots\}$, such that
$\{e_1,e_2,\dots\}$, such that
	\[
	e_i \cdot e_j=
	\begin{cases}
	0, &\text{if \ $i\ne j$;}\\
	\displaystyle \sum_{k}a_{ik}e_k, &\text{if \ $i=j$,}
	\end{cases}
	\]
then this algebra is called an evolution algebra. The basis is called a natural basis.
\end{defn}

The matrix $A=(a_{ij})$ is called the matrix of structural constants.

Evolution algebras have the following elementary properties (see \cite{t}).
Evolution algebras are not
associative, in general; they are commutative, flexible, but not power-associative, in general; direct
sums of evolution algebras are also evolution algebras; Kronecker products of evolutions algebras are
also evolution algebras.
The next theorem gives the classification of the real
two-dimensional evolution algebras.

\begin{thm}[\cite{uzmj}]\label{ms}
 Any two-dimensional real evolution algebra \textit{E} is isomorphic
	to one of the following pairwise non-isomorphic algebras:
\begin{itemize}
\item[(i)]  $\dim E^2=1$.\\
	$E_{1}: e_1e_1=e_1, \ \ e_2e_2=0$;\\
	$E_{2}: e_1e_1=e_1, \ \  e_2e_2=e_1$;\\
	$E_{3}: e_1e_1=e_1+e_2, \ \  e_2e_2=-e_1-e_2$;\\
	$E_4: e_1e_1=e_2, \ \  e_2e_2=0 $;\\
	$E_5: e_1e_1=e_2, \ \  e_2e_2=-e_2$;
	
\item[(ii)] $\dim E^2=2$.\\
$E_6(a_2;a_3): e_1e_1=e_1+a_2 e_2, \ \ e_2e_{2}=a_3 e_1+e_2, \ \ 1-a_2 a_3 \neq 0, a_2, a_3\in {\mathbb{R}}$. Moreover $E_6(a_2;a_3)$
	is isomorphic to $E_6(a_3;a_2)$.\\
	$E_7(a_4): e_1e_1=e_2, \ \ e_2e_{2}=e_1+a_4 e_2, \ \ where \ a_4\in {\mathbb{R}}$.
\end{itemize}
\end{thm}

Following \cite{clr} we consider a family $\left\{E^{[s,t]}:\ s,t \in \mathbb{R},\ 0\leq s\leq t
\right\}$ of $n$-dimensional evolution algebras over the field $\mathbb{R}$,
with basis $e_1,\dots,e_n$, and the multiplication table
\[
e_ie_i =\sum_{j=1}^na_{ij}^{[s,t]}e_j, \ \ i=1,\dots,n; \qquad
e_ie_j =0,\ \ i\ne j.
\]
 Here parameters $s,t$ are
considered as time.

Denote by
$M^{[s,t]}=\left(a_{ij}^{[s,t]}\right)_{i,j=1,\dots,n}$ the matrix
of structural constants.

\begin{defn}
 A family $\left\{E^{[s,t]}:\ s,t \in \mathbb{R},\ 0\leq s\leq t
	\right\}$ of $n$-dimensional evolution algebras over the field $\mathbb{R}$
	is called a chain of evolution algebras (CEA) if the matrix
	$M^{[s,t]}$ of structural constants satisfies the
	Chapman-Kolmogorov equation
\begin{equation}\label{2}
	M^{[s,t]}=M^{[s,\tau]}M^{[\tau,t]}, \ \ \text{for any} \ \
	s<\tau<t.
\end{equation}
\end{defn}

\begin{defn}%\label{d4}
Assume a CEA, $E^{[s,t]}$, has a property, say $P$,
	at pair of times $(s_0,t_0)$; one says that the CEA has $P$ property
	transition if there is a pair $(s,t)\ne (s_0,t_0)$ at which the CEA
	has no the property $P$.
\end{defn}

Denote
\begin{align*}
\mathcal T & =\{(s,t): 0\leq s\leq t\}; \\
\mathcal T_P & =\{(s,t)\in \mathcal T: E^{[s,t]} \ \ \text{has  property} \  P \}; \\
\mathcal T_P^0 & =\mathcal T\setminus \mathcal T_P=\{(s,t)\in \mathcal T: E^{[s,t]} \ \ \text{has no  property} \ P \}.
\end{align*}
The sets have the following meaning:
\begin{itemize}
\item $\mathcal T_P$-the duration of the property $P$;

\item  $\mathcal T_P^0$-the lost duration of the property $P$;
\end{itemize}

The partition $\{\mathcal T_P, \mathcal T^0_P\}$ of the set
$\mathcal T$ is called $P$ property diagram.

For example, if $P=$ commutativity then since any evolution algebra
is commutative, we conclude that any CEA has not commutativity
property transition.

To construct a chain of two-dimensional evolution algebras one has to solve equation \eqref{2}
for the $2\times 2$ matrix $\M^{[s,t]}$. This equation gives the following system of functional equations (with four unknown functions):
\begin{equation}\label{fe}
\begin{array}{llll}
a_{11}^{[s,t]}=a_{11}^{[s,\tau]}a_{11}^{[\tau,t]}+a_{12}^{[s,\tau]}a_{21}^{[\tau,t]},\\[2mm]
a_{12}^{[s,t]}=a_{11}^{[s,\tau]}a_{12}^{[\tau,t]}+a_{12}^{[s,\tau]}a_{22}^{[\tau,t]},\\[2mm]
a_{21}^{[s,t]}=a_{21}^{[s,\tau]}a_{11}^{[\tau,t]}+a_{22}^{[s,\tau]}a_{21}^{[\tau,t]},\\[2mm]
a_{22}^{[s,t]}=a_{21}^{[s,\tau]}a_{12}^{[\tau,t]}+a_{22}^{[s,\tau]}a_{22}^{[\tau,t]}.\\
 \end{array}
 \end{equation}

But the general analysis of system \eqref{fe} is difficult.

In \cite{mu} we studied  the classification dynamics of known two-dimensional chains of evolution algebras constructed in \cite{rm} and showed that known chains of evolution algebras never contain evolution algebra isomorphic to $E_4$ in any time $s, t$ (see Theorem~\ref{ms}). In this section we will construct CEAs, which will contain the algebra $E_4$ for some period of time.

To construct a CEA, which for some time will be isomorphic to $E_4$, we need the following theorem.

\begin{thm}[\cite{mu}]
	An evolution algebra $E_{\mathcal{M}}$ is isomorphic to $E_4$ if and only if
	$E_{\mathcal{M}}$ has the matrix of structural constants in the following form:
\begin{equation}\label{MM}
	\M_1=\begin{pmatrix}
	0 & \beta\\
	0 & 0
	\end{pmatrix} \ \ \text{or} \ \
	\M_2=\begin{pmatrix}
	0 & 0\\
	\gamma & 0
	\end{pmatrix}, \quad  \text{where} \  \ \beta,\gamma \in \mathbb{R}.
\end{equation}
\end{thm}

Thus, we should construct CEAs with the matrix of structural constant showed in \eqref{MM}.

Consider \eqref{fe} with $a_{11}^{[s,t]}=\alpha(s,t), \ a_{12}^{[s,t]}=\beta(s,t), \ a_{21}^{[s,t]}=\gamma(s,t), \ a_{22}^{[s,t]}=\delta(s,t)$. Therefore, to find CEA, we should solve the next equation:

\begin{equation}\label{c11}
\begin{pmatrix}
{\alpha(s,\tau)} & {\beta(s,\tau)}\\[2mm]
{\gamma(s,\tau)} & {\delta(s,\tau)}
\end{pmatrix}
\cdot
\begin{pmatrix}
{\alpha(\tau,t)} & {\beta(\tau,t)}\\[2mm]
{\gamma(\tau,t)} & {\delta(\tau,t)}
\end{pmatrix}
=
\begin{pmatrix}
{\alpha(s,t)} & {\beta(s,t)}\\[2mm]
{\gamma(s,t)} & {\delta(s,t)}
\end{pmatrix}.
\end{equation}

\textbf{Case 1.1.} If we consider in \eqref{c11}, $\alpha(s, t)= \gamma(s, t)\equiv 0, \beta(s,t)\neq 0, \delta(s,t)\neq 0$, then we have the following:

\begin{equation}\label{c1}
\begin{pmatrix}
{0} & {\beta(s,\tau)}\\[2mm]
{0} & {\delta(s,\tau)}
\end{pmatrix}
\cdot
\begin{pmatrix}
{0} & {\beta(\tau,t)}\\[2mm]
{0} & {\delta(\tau,t)}
\end{pmatrix}
=
\begin{pmatrix}
{0} & {\beta(s,t)}\\[2mm]
{0} & {\delta(s,t)}
\end{pmatrix}.
\end{equation}

From \eqref{c1}, we get the following system of functional equations:

\begin{equation}\label{s11}
\left\{\begin{array}{llllll}
\beta(s,\tau)\delta(\tau,t)=\beta(s,t), \\[3mm]
\delta(s,\tau)\delta(\tau,t)=\delta(s,t).
\end{array}
\right.
\end{equation}

The second equation of system \eqref{s11} is known as Cantor's second equation, which has the following solutions:
\begin{itemize}
  \item[(1)] $\delta(s,t)\equiv 0$;
   \item[(2)] $\delta(s,t)=\frac{\phi(t)}{\phi(s)}$, where $\phi$ is an arbitrary function with $\phi(s)\neq 0$;
   \item[(3)] $\delta(s,t)= \begin{cases}
1, & \text{if} \ \ 0<s\leq t<a; \\
0, &\text{if} \ \ t\geq a.
\end{cases}$
\end{itemize}

Substituting these solutions in the first equation of \eqref{s11}, we find $\beta(s,t)$:

\begin{itemize}
  \item[(1)] $\beta(s,t)\equiv 0$;
   \item[(2)] $\beta(s,t)=\rho(s)\phi(t)$, where $\rho$ is an arbitrary function;
   \item[(3)] $\beta(s,t)= \begin{cases}
\sigma(s), & \text{if} \ \ 0<s\leq t<a; \\
0, & \text{if} \ \ t\geq a,
\end{cases}$
\end{itemize}
where $\sigma$ is an arbitrary function;

From these solutions we have the following matrices of structural constants of CEAs:

\begin{align*}
\M_0^{[s,t]}& =\begin{pmatrix}
{0} & {0}\\[2mm]
{0} & {0}
\end{pmatrix},\\
\M_1^{[s,t]}&=\begin{pmatrix}
{0} & {\rho(s)\phi(t)}\\[2mm]
{0} & {\frac{\phi(t)}{\phi(s)}}
\end{pmatrix},
\end{align*}
 where  $\rho, \phi$ are arbitrary functions, with $\phi(s)\neq0$;

\[\M_2^{[s,t]}=\begin{cases}
\begin{pmatrix}
0 & \sigma(s)\\
0 & 1
\end{pmatrix}, &  \text{if} \ \ 0<s\leq t<a;\\[4mm]
\begin{pmatrix}
{0} & {0}\\
{0} & {0}
\end{pmatrix},& \text{if} \ \ t\geq a,\\
\end{cases}\]
 where $a>0$ and $\sigma$ is an arbitrary function.

\textbf{Case 1.2.} Consider the case $\alpha(s, t)= \beta(s, t)\equiv 0, \gamma(s,t)\neq 0, \delta(s,t)\neq 0$. Then from \eqref{c11} we have the following:

\begin{equation*} % \label{c12}
\begin{pmatrix}
{0} & {0}\\[2mm]
{\gamma(s,\tau)} & {\delta(s,\tau)}
\end{pmatrix}
\cdot
\begin{pmatrix}
{0} & {0}\\[2mm]
{\gamma(\tau,t)} & {\delta(\tau,t)}
\end{pmatrix}
=
\begin{pmatrix}
{0} & {0}\\[2mm]
{\gamma(s,t)} & {\delta(s,t)}
\end{pmatrix}.
\end{equation*}

From the last equality, we have the following system of equations:

\begin{equation}\label{s12}
\left\{\begin{array}{llllll}
\delta(s,\tau)\gamma(\tau,t)=\gamma(s,t), \\[3mm]
\delta(s,\tau)\delta(\tau,t)=\delta(s,t).
\end{array}
\right.
\end{equation}

The second equation of  system \eqref{s12}, known as Cantor's second equation, which has the following solutions:

(1) $\delta(s,t)\equiv 0$;

(2) $\delta(s,t)=\frac{\varphi(t)}{\varphi(s)}$, where $\varphi$ is an arbitrary function with $\varphi(s)\neq 0$;

(3)  $\delta(s,t)=\begin{cases}
1, & \text{if} \ \ 0<s\leq t<a; \\
0, & \text{if} \ \ t\geq a.
\end{cases}
$

Substituting these solutions in the first equation of \eqref{s12}, we find $b(s,t)$:

(1) $\gamma(s,t)\equiv 0$;

(2) $\gamma(s,t)=\frac{f(t)}{\varphi(s)}$, where $f$ is an arbitrary function;

(3)  $\gamma(s,t)=\begin{cases}
g(t), & \text{if} \ \ 0<s\leq t<a; \\
0, & \text{if} \ \ t\geq a.
\end{cases}$ \quad where $g$ is an arbitrary function.

From these solutions we have the next matrices of structural constants of CEAs:

\[\M_0^{[s,t]}=\begin{pmatrix}
{0} & {0}\\[2mm]
{0} & {0}
\end{pmatrix},\]

\[
\M_3^{[s,t]}=\begin{pmatrix}
{0} & {0}\\[2mm]
{\frac{f(t)}{\varphi(s)}} & {\frac{\varphi(t)}{\varphi(s)}}
\end{pmatrix},\]
 where  $f, \varphi$ are arbitrary functions, $\varphi(s)\neq0$;

\[\M_4^{[s,t]}=\begin{cases}
\begin{pmatrix}
0 & 0\\
g(t) & 1
\end{pmatrix}, & \text{if} \ \ 0<s\leq t<a;\\[4mm]
\begin{pmatrix}
{0} & {0}\\
{0} & {0}
\end{pmatrix}, & \text{if} \ \ t\geq a,\\
\end{cases}\] where $a>0$ and $g$ is an arbitrary function.

\textbf{Case 1.3.} Let us try to find the solution satisfying the following:

\begin{equation}\label{c13}
\begin{pmatrix}
{\alpha(s,\tau)} & {\beta(s,\tau)}\\[2mm]
{\gamma(s,\tau)} & {\delta(s,\tau)}
\end{pmatrix}
\cdot
\begin{pmatrix}
{\alpha(\tau,t)} & {\beta(\tau,t)}\\[2mm]
{\gamma(\tau,t)} & {\delta(\tau,t)}
\end{pmatrix}
=
\begin{pmatrix}
{0} & {\beta(s,t)}\\[2mm]
{0} & {0}
\end{pmatrix}.
\end{equation}

From \eqref{c13} we have the next system of functional equations:

\begin{equation}\label{sec13}
\left\{\begin{array}{llllll}
\alpha(s,\tau)\alpha(\tau,t)+\beta(s,\tau)\gamma(\tau,t)=0, \\[3mm]
\alpha(s,\tau)\beta(\tau,t)+\beta(s,\tau)\delta(\tau,t)=\beta(s,t), \\[3mm]
\gamma(s,\tau)\alpha(\tau,t)+\delta(s,\tau)\gamma(\tau,t)=0, \\[3mm]
\gamma(s,\tau)\beta(\tau,t)+\delta(s,\tau)\delta(\tau,t)=0.
\end{array}
\right.
\end{equation}

Let $\alpha(s, t)=\gamma(s,t)=0$. Then we get:

\begin{equation}\label{sc13}
\left\{\begin{array}{llllll}
\beta(s,\tau)\delta(\tau,t)=\beta(s,t), \\[3mm]
\delta(s,\tau)\delta(\tau,t)=0.
\end{array}
\right.
\end{equation}

To find a non-zero solution of the system of equations \eqref{sc13} we should prove that
the equation
\begin{equation}\label{ds12}
\delta(s,\tau)\delta(\tau,t)=0, \ \ \ \text{for all} \ \  s<\tau<t,
\end{equation}
has a non-zero solution. Indeed, take $C>0$ and
\begin{equation}\label{f1}
\delta(s,t)=\begin{cases}
0, & \text{if} \ \ 0<C\leq s<t \ \ \text{or} \ \  0<s<t\leq C;\\[2mm]
f(s,t), & \text{if} \ \  0<s<C<t,
\end{cases}
\end{equation}
where $f(s,t)$ is an arbitrary non-zero function.

Now, we show that independently on $f(s,t)$ the function \eqref{f1} satisfies \eqref{ds12}:
take an arbitrary $\tau$ such that $s<\tau<t$, then for given $C>0$,
we have only two possibilities:

\textbf{Case 1.3.1.} Let $\tau\leq C$. By the defined function \eqref{f1}, we have that $\delta(s,\tau)=0$
and for $\delta(\tau, t)$:

\begin{equation}\label{f112}
\delta(\tau,t)=\begin{cases}
0, & \text{if} \  \ t\leq C;\\
f(\tau,t), &\text{if} \ \  t>C,
\end{cases}
\end{equation}
where $f(\tau,t)$ is the function fixed in \eqref{f1}.

Therefore, $\delta(s,\tau)\delta(\tau,t)=0$.

\textbf{Case 1.3.2.} $\tau>C$. Also from \eqref{f1}, we have that $\delta(\tau,t)=0$
and for $\delta(s,\tau)$:

\begin{equation*}%\label{f12}
\delta(s,\tau)=\begin{cases}
f(s,\tau),& \text{if} \ \  s<C;\\
0, & \text{if} \ \ s\geq C,
\end{cases}
\end{equation*}
where $f(s, \tau)$ is the function fixed in \eqref{f1}.

Therefore, $\delta(s,\tau)\delta(\tau,t)=0$.

Thus, we have proved that the function \eqref{f1} satisfies  equation \eqref{ds12}.

Now we should find solutions of the first equation of system \eqref{sc13}:

\begin{equation}\label{bs13}
\beta(s,\tau)\delta(\tau,t)=\beta(s,t), \ \  s<\tau<t,
\end{equation}
where $\delta(\tau,t)$ is given by \eqref{f1}.

To find a solution we have the next possibilities:

 \textbf{Case 1.3.3.} Let $\tau \leq C$, then by the defined function \eqref{f1} we have that $\delta(s,\tau)=0$ and from \eqref{f112} in period of time $t\leq C$, $\delta(\tau, t)=0$, then from \eqref{bs13} we have that $\beta(s,t)=0$. When $t>C$, $\delta(\tau,t)=f(\tau,t)$ and by \eqref{bs13} we have to solve the next equation:

 \begin{equation}\label{bs131}
 \beta(s,\tau)f(\tau,t)=\beta(s,t), \ \  s<\tau<t.
 \end{equation}

 We solve \eqref{bs131} for some particular cases:

 \textbf{Case 1.3.3.1} Consider $\beta(s,t)=f(s,t)$. Then from \eqref{bs131}, we have $f(s,\tau)f(\tau,t)=f(s,t)$, which is Cantor's second equation.
 As $f(s,t)$ is a non-zero function, then we have the next solution:

 \[f(s,t)=\frac{\Phi(t)}{\Phi(s)},\]
 where $\Phi$ is an arbitrary function, with $\Phi(s) \neq 0$.

 Thus we have the next solution of system \eqref{sec13}:

 \begin{equation*} %\label{y111}
 \begin{array}{llll}
 \alpha(s, t)\equiv 0,\\
 \beta(s, t)=\begin{cases}
 0, &  \text{if} \ \ s<t \leq C\\
 \frac{\Phi(t)}{\Phi(s)}, & \text{if} \ \  t>C,
\end{cases}\\
 \gamma(s, t)\equiv 0,\\
 \delta(s, t)=\begin{cases}
 0, & \text{if} \ \ 0<C\leq s< t \ \ \text{or} \ \ 0<s<t\leq C; \\
 \frac{\Phi(t)}{\Phi(s)}, & \text{if} \ \ s< C<t,\\
 \end{cases}
 \end{array}
 \end{equation*}
where $C>0$ and $\Phi$ is an arbitrary function, with $\Phi(s)\neq 0$.

 Then we have the next matrix of structural constants:

\[\M_5^{[s,t]}=\begin{cases}
\begin{pmatrix}
0 & 0\\
0 & 0
\end{pmatrix}, & \text{if} \ \ \ \ s<t \leq C;\\[4mm]
\begin{pmatrix}
0 & \frac{\Phi(t)}{\Phi(s)}\\
0 & 0
\end{pmatrix}, & \text{if} \ \ \ t>C,
\end{cases}\] where $C>0$ and $\Phi$ is an arbitrary function, with $\Phi(t)\neq 0$.

 \textbf{Case 1.3.3.2.} Let $\beta(s,t)\neq f(s,t)$.  As $f(\tau,t)$ is an arbitrary non-zero function,
 consider $f(\tau,t)=\frac{\phi(\tau)}{\phi(t)}$, with $\phi(t)\neq 0$. Then from \eqref{bs131} we have the following:

 \begin{align*}
 \beta(s,\tau)\cdot \frac{\phi(\tau)}{\phi(t)}& =\beta(s,t),\\
\beta(s,t)\phi(t) &=\beta(s,\tau)\phi(\tau).
\end{align*}

From the last equality, we can see that $\beta(s,t)\phi(t)$ does not depend on $t$,
i.e. there exists a function $\rho(s)$ such that $\beta(s,t)\phi(t)=\rho(s)$. Therefore, $\beta(s,t)=\frac{\rho(s)}{\phi(t)}$.

Then we get the next solution of  system \eqref{sec13}:

\begin{equation*} %\label{y11}
\begin{array}{llll}
\alpha(s, t)\equiv 0,\\
\beta(s, t)=\begin{cases}
0, & \text{if} \ \ s<t \leq C;\\
\frac{\rho(s)}{\phi(t)}, & \text{if} \ \  t>C,
\end{cases}\\
\gamma(s, t)\equiv 0,\\
\delta(s, t)=\begin{cases}
0, & \text{if} \ \ 0<C\leq s< t \ \ \text{or} \ \ 0<s<t\leq C; \\
\frac{\phi(s)}{\phi(t)}, & \text{if} \ \ s< C<t,\\
\end{cases}
\end{array}
\end{equation*}
where $C>0$ and $\phi, \rho$ are arbitrary functions with $\phi(t)\neq 0$.

Then we have, respectively, the next matrix of structural constants  to the solution:

\[\M_6^{[s,t]}=\begin{cases}
\begin{pmatrix}
0 & 0\\
0 & 0
\end{pmatrix}, & \text{if} \ \  \ s<t \leq C;\\[4mm]
\begin{pmatrix}
0 & \frac{\rho(s)}{\phi(t)}\\
0 & 0
\end{pmatrix}, &\text{if} \ \ \ t>C,
\end{cases}\]
where $C>0$ and $\phi, \rho$ are arbitrary functions with $\phi(t)\neq 0$.

\textbf{Case 1.3.4.} When $\tau > C$, then by the defined function \eqref{f1} we have that $\delta(\tau, t)=0$. Then from \eqref{bs13}
we have that $\beta(s, t)=0$. Thus we get the trivial CEA.

\textbf{Case 1.4.} Let us try to find the solution satisfying:

\begin{equation}\label{c2}
\begin{pmatrix}
{\alpha(s,\tau)} & {\beta(s,\tau)}\\[2mm]
{\gamma(s,\tau)} & {\delta(s,\tau)}
\end{pmatrix}
\cdot
\begin{pmatrix}
{\alpha(\tau,t)} & {\beta(\tau,t)}\\[2mm]
{\gamma(\tau,t)} & {\delta(\tau,t)}
\end{pmatrix}
=
\begin{pmatrix}
{0} & {0}\\[2mm]
{\gamma(s,t)} & {0}
\end{pmatrix}.
\end{equation}
From equality \eqref{c2} we have the next system of functional equations:

\begin{equation*}%\label{sc2}
\left\{\begin{array}{llllll}
\alpha(s,\tau)\alpha(\tau,t)+\beta(s,\tau)\gamma(\tau,t)=0, \\[3mm]
\alpha(s,\tau)\beta(\tau,t)+\beta(s,\tau)\delta(\tau,t)=0, \\[3mm]
\gamma(s,\tau)\alpha(\tau,t)+\delta(s,\tau)\gamma(\tau,t)=\gamma(s, t), \\[3mm]
\gamma(s,\tau)\beta(\tau,t)+\delta(s,\tau)\delta(\tau,t)=0.
\end{array}
\right.
\end{equation*}

Let $\alpha(s,t)=\beta(s,t)=0$. Then we have the next system:

\begin{equation*}%\label{sc14}
\left\{\begin{array}{llllll}
\delta(s,\tau)\gamma(\tau,t)=\gamma(s,t), \\[3mm]
\delta(s,\tau)\delta(\tau,t)=0.
\end{array}
\right.
\end{equation*}

The analysis of this system is similar to the \eqref{sc13} and we get the following CEAs:

\[\M_7^{[s,t]}=\begin{cases}
\begin{pmatrix}
0 & 0\\
\frac{\Psi(t)}{\Psi(s)} & 0
\end{pmatrix}, & \text{if} \ \ \  s<C;\\[4mm]
\begin{pmatrix}
0 & 0\\
0 & 0
\end{pmatrix}, & \text{if} \ \ \, s\geq C,
\end{cases}\] where $C>0$ and $\Psi$ is an arbitrary function, with $\Psi(t)\neq 0$;

\[\M_8^{[s,t]}= \begin{cases}
\begin{pmatrix}
0 & 0\\
\frac{\sigma(t)}{\varphi(s)} & 0
\end{pmatrix}, & \text{if} \ \ \ s<C;\\[4mm]
\begin{pmatrix}
0 & 0\\
0 & 0
\end{pmatrix}, & \text{if} \ \ \ s\geq C,
\end{cases}\] where $C>0$ and $\varphi, \sigma$ are arbitrary functions with $\varphi(s)\neq 0$.

Denote by $E_i^{[s,t]}$ the CEA with matrix $\M_i^{[s,t]}$.

\begin{rk} We should note that, from the CEAs $E_i^{[s,t]}, i=1,\dots,8$, only $E_3^{[s,t]}$ coincides
with the CEA $E_{16}^{[s,t]}$ constructed in \cite{rm} and it has the same dynamic. All other CEAs are different
from CEAs constructed in \cite{rm} and have different dynamics.
\end{rk}

The next theorem gives time-depending dynamics of these CEAs:	

\begin{thm} For the next CEAs  hold:
	
\[E_{1}^{[s,t]} \simeq \begin{cases}
E_1 & \text{for all} \ \ (s,t)\in \left\{(s,t): s<t, \ \ \rho(s)=0 \right\}, \\
E_2 & \text{for all} \ \ (s,t)\in \left\{(s,t): s<t, \ \ \rho(s)\neq 0 \right\}.
\end{cases}\]
		
\[E_{2}^{[s,t]} \simeq \begin{cases}
E_1 & \text{for all} \ \ (s,t)\in \left\{(s,t): s<t<a, \ \ \sigma(s)=0 \right\}, \\
E_2 & \text{for all} \ \ (s,t)\in \left\{(s,t): s<t<a, \ \ \sigma(s)\neq 0 \right\},\\
E_0 &\text{for all} \ \ (s,t)\in \left\{(s,t): t\geq a\right\}
\end{cases}\]

$E_{3}^{[s,t]}$ is isomorphic to $E_1$ for any $(s,t)\in \mathcal T$.

\[E_{4}^{[s,t]} \simeq \begin{cases}
E_1 &  \text{for all} \ \ (s,t)\in \left\{(s,t): s<t<a\right\},\\
E_0 &  \text{for all} \ \ (s,t)\in \left\{(s,t): t \geq a\right\}.
\end{cases}\]

\[E_{5}^{[s,t]} \simeq \begin{cases}
E_0 &  \text{for all} \ \ (s,t)\in \left\{(s,t): s<t\leq C \right\},\\
E_4 &  \text{for all} \ \ (s,t)\in \left\{(s,t): t > C\right\}.
\end{cases}\]

\[E_{6}^{[s,t]} \simeq \begin{cases}
E_0 &  \text{for all} \ \ (s,t)\in \left\{(s,t): s<t\leq C \right\},\\
E_0 &  \text{for all} \ \ (s,t)\in \left\{(s,t): t> C, \ \ \rho(s)=0\right\},\\
E_4 &  \text{for all} \ \ (s,t)\in \left\{(s,t): t>C, \ \ \rho(s)\neq 0 \right\}.
\end{cases}\]

\[E_{7}^{[s,t]} \simeq \begin{cases}
E_4 &  \text{for all} \ \ (s,t)\in \left\{(s,t): s<C\right\},\\
E_0 &  \text{for all} \ \ (s,t)\in \left\{(s,t): s\geq C\right\}.
\end{cases}\]

\[E_{8}^{[s,t]} \simeq \begin{cases}
E_0 & \text{for all} \ \ (s,t)\in \left\{(s,t): s<C, \ \ \sigma(t)=0 \right\},\\
E_4 &  \text{for all} \ \ (s,t)\in \left\{(s,t): s<C, \ \ \sigma(t)\neq 0 \right\}, \\
E_0 &  \text{for all} \ \ (s,t)\in \left\{(s,t): s\geq C \right\}.
\end{cases}\]
	\end{thm}

\proof
When $\rho(s)=0$, then $E_1^{[s,t]}$ is isomorphic to $E_1$,  for all $s, t \in \mathcal T$ by the change of basis $e'_1=e_1,  \ e'_2=\frac{\phi(s)}{\phi(t)}e_2$, and when $\rho(s)\neq 0$, it is isomorphic to
$E_2$,  for all $s, t \in \mathcal T$ by the change of basis $e_1'=\frac{1}{\rho(s)\phi(t)}e_1, \ e'_2=\frac{\phi(s)}{\phi(t)}e_2$.

When $\sigma(s)=0$, then $E_2^{[s,t]}$ is isomorphic to $E_1$, for all $s, t \in \mathcal T, \ s<t<a$, by the change of basis  $e'_1=e_1,  \ e'_2=e_2$, and when $\sigma(s)\neq 0$, it is isomorphic to
$E_2$,  for all $s, t \in \mathcal T, \ s<t<a$, by the change of basis $e_1'=\frac{1}{\sigma(s)}e_1, \ e'_2=e_2$. In period of time $t\geq a$, it  will
be isomorphic to the trivial evolution algebra $E_0$.

$E_3^{[s,t]}$  is isomorphic to $E_1$, for all $s, t \in \mathcal T$ by the change of basis
$e'_2=\frac{f(t)\varphi(s)}{\varphi^2(t)}e_1+\frac{\varphi(s)}{\varphi(t)}e_2, \ e'_2=e_1$.

$E_4^{[s,t]}$ is isomorphic to $E_1$, for all $s, t \in \mathcal T, \ s<t<a$, by the change of basis  $e'_1=\sigma(t)e_1+e_2, \ \ e'_2=e_1$, in period of time $t\geq a$, it  will be isomorphic to the trivial evolution algebra $E_0$.

$E_5^{[s,t]}$ is isomorphic to $E_4$,  for all $s, t \in \mathcal T, \ t>C$, by the change of basis  $e'_1=\frac{\Phi(s)}{\Phi(t)}e_1, \ \ e'_2=e_2$, in period of time $s<t\leq C$, it  will be isomorphic to the trivial evolution algebra $E_0$.

When $\rho(s)\neq 0$, then $E_6^{[s,t]}$ is isomorphic to $E_4$, for all $s, t \in \mathcal T, \ t>C$, by the change of basis  $e'_1=\frac{\phi(t)}{\rho(s)}e_1, \ \ e'_2=e_2$, in period of time $s<t\leq C$ and when $\rho(s)=0$, then it will be isomorphic to the trivial evolution algebra $E_0$.

$E_7^{[s,t]}$ is isomorphic to $E_4$, for all $s, t \in \mathcal T, \ s<C$, by  the change of basis  $e'_1=\frac{\Psi(s)}{\Psi(t)}e_1, \ \ e'_2=e_2$, in period of time $s\geq C$, it  will be isomorphic to the trivial evolution algebra $E_0$.

When $\sigma(t)\neq 0$, then $E_6^{[s,t]}$ is isomorphic to $E_4$, for all $s, t \in \mathcal T, \ s<C$, by the change of basis  $e'_1=\frac{\varphi(s)}{\sigma(t)}e_1, \ \ e'_2=e_2$, in period of time $s\geq C$, and when $\sigma(t)=0$, then it will be isomorphic to the trivial evolution algebra $E_0$.
\endproof

Thus we proved that there exists CEAs that for some values of time will be isomorphic to $E_4$.

\section{Rota-Baxter operators on Evolution Algebras}\label{S:rota}

In this section we will study Rota-Baxter operators on evolution algebras.

\begin{defn}
Let $\F$ be a field. A Rota-Baxter operator of weight $\lambda \in \F$ on evolution algebra $(E, \cdot)$ over $\F$ is a linear map
$P:E\longrightarrow E$ satisfying

\[
P(x)\cdot P(y)=P(x\cdot P(y)+P(x)\cdot y+\lambda\cdot x\cdot y), \ \ \text{for all} \ \ x,y \in E.
\]

\end{defn}

Note that, if $P$ is a Rota-Baxter operator of weight $\lambda \neq 0$, then $\lambda^{-1}P$ is a Rota-Baxter
operator $P$ of weight 1. Therefore, one only needs to consider Rota-Baxter operators of
weight 0 and 1. We also assume that $\F=\C$.

To study Rota-Baxter operators on an evolution algebra $(E, \cdot)$ over the field $\C$, we need to the next theorem, which gives the classification
of two-dimensional complex evolution algebras.

\begin{thm}[\cite{clor}]
 Any two-dimensional complex evolution algebra $E$ is isomorphic
	to one of the following pairwise non-isomorphic algebras:
	
	(i) $\dim E^2=1$
	
	\begin{itemize}
		\item 	$E_{1}: e_1e_1=e_1$;
		\item   $E_{2}: e_1e_1=e_1, \ \  e_2e_2=e_1$;
		\item   $E_{3}: e_1e_1=e_1+e_2, \ \  e_2e_2=-e_1-e_2$;
		\item	$E_4: e_1e_1=e_2$;
	\end{itemize}
	
	(ii) $\dim E^2=2$:
	\begin{itemize}
		\item 	$E_5(a_2,a_3): e_1e_1=e_1+a_2 e_2, \ \ e_2e_{2}=a_3 e_1+e_2, \ \ 1-a_2 a_3 \neq 0, \ a_2, a_3\in {\mathbb{C}}$. Moreover $E_5(a_2,a_3)\cong E_5(a_3,a_2)$.
		\item 	$E_6(a_4): e_1e_1=e_2, \ \ e_2e_{2}=e_1+a_4 e_2$, where $E_6(a_4)\cong E_6(a_4') \Leftrightarrow \frac{a_4'}{a_4}=\cos \frac{2\pi k}{3}+i \sin \frac{2 \pi k}{3}$ for some $k=0,1,2$.
	\end{itemize}
\end{thm}

For an $n$-dimensional evolution algebra we have that $e_i e_i=\sum_{j=1}^{n}a_{ij}e_j$,  for all $i$,
$e_i e_j=0, \, i \neq j$, and for the Rota-Baxter operator  $P(e_i)=\sum_{j=1}^{n}r_{ij}e_j$.

To find a Rota-Baxter operator of weight  $\lambda$ on evolution algebra we should solve the system of equations obtained from the next equality:

\begin{equation}\label{rb}
P(e_i) P(e_j) = P(e_i P(e_j)+P(e_j) e_i+\lambda e_i e_j), \ \ i, j=1,\dots,n.
\end{equation}

For $i=j$, the LHS of \eqref{rb} equals to:

\begin{equation}\label{rb1}
P(e_i) P(e_i)=\sum_{j=1}^{n} \sum_{k=1}^{n}r_{ij}^2 a_{jk}e_k.
\end{equation}

And the RHS \eqref{rb}, for $i=j$, is equal to:

\[
P(e_i P(e_i)+P(e_i) e_i+\lambda e_i e_i)=(2r_{ii}+\lambda)\sum_{j=1}^{n}\sum_{k=1}^{n}a_{ij}r_{jk}e_k.
\]

For $i\neq j$, we have:

\[
P(e_i) P(e_j)=\sum_{k=1}^{n} \sum_{s=1}^{n}r_{ik} r_{jk} a_{ks}e_s
\]
\begin{equation}\label{rb22}
P(e_i P(e_j)+P(e_i) e_j+\lambda e_i e_j)=r_{ji}\sum_{k=1}^{n}\sum_{l=1}^{n}a_{ik}r_{kl}e_l+r_{ij}\sum_{k=1}^{n}\sum_{l=1}^{n}a_{js}r_{sl}e_l
\end{equation}

Thus, to find an operator of Rota-Baxter on an evolution algebra we have to solve  the system of equations \eqref{rb} with equalities \eqref{rb1}--\eqref{rb22}.

Without loss of generality for the Rota-Baxter operator with $\lambda \neq 0$ we can assume $\lambda=1$, therefore we will find the Rota-Baxter operators corresponding to a given evolution algebra with weights $\lambda=0$ and $\lambda=1$.

Now, we will find  the Rota-Baxter operators corresponding to the two-dimensional complex evolution algebras $E_i, \ i=1,\dots, 6$, with weights $\lambda=0$ and $\lambda=1$.

Consider the Rota-Baxter operator
\[\begin{pmatrix}
P(e_1)\\[2mm]
P(e_2)
\end{pmatrix}
=
\begin{pmatrix}
a & b\\[2mm]
c & d
\end{pmatrix}
\begin{pmatrix}
e_1\\[2mm]
e_2
\end{pmatrix},
\]
where $a, b, c, d \in \mathbb{C}$.

\begin{itemize}
	\item For the algebra $E_{1}: e_1e_1=e_1$ to find the matrix form of the Rota-Baxter operator of weight $0$ we should solve the next system of equations:
	
	\[\left\{\begin{array}{llllll}
	a^2=0, \\[1mm]
	2ab=0, \\[1mm]
	c^2=0, \\[1mm]
	ac=ac, \\[1mm]
	bc=0.
	\end{array}
	\right.\]
	
	This system of equations has the following solution:	
	
	\[\begin{pmatrix}
	0 & b\\[2mm]
	0 & d
	\end{pmatrix}.\]
	
	For the algebra $E_{1}: e_1e_1=e_1$ to find the matrix form of the Rota-Baxter operator of weight $\lambda=1$ we should solve the next system of equations:
	
	\[\left\{\begin{array}{llllll}
	a^2+a=0, \\[1mm]
	(2a+1)b=0, \\[1mm]
	c^2=0, \\[1mm]
	ac=ac, \\[1mm]
	bc=0.
	\end{array}
	\right.\]
	
	This system of equations has the following solutions:	
	
	\[\begin{pmatrix}
	-1 & 0\\[2mm]
	0 & d
	\end{pmatrix}, \qquad
	\begin{pmatrix}
	0 & 0\\[2mm]
	0 & d
	\end{pmatrix}.\]
	
	\item For the algebra $E_{2}: e_1e_1=e_1, \ \  e_2e_2=e_1$ to find the matrix form of the Rota-Baxter operator of weight $0$ we should solve the next system of equations:
	
	\[\left\{\begin{array}{llllll}
	b^2=a^2, \\[1mm]
	2ab=0, \\[1mm]
	c^2+d^2=2ad, \\[1mm]
	2bd=0, \\[1mm]
	bd=ab, \\[1mm]
	(b+c)b=0.
	\end{array}
	\right.\]
	
	This system of equations has the next solutions:
	
	\[\begin{pmatrix}
	0 & 0\\[2mm]
	c & -ic
	\end{pmatrix}, \qquad
	\begin{pmatrix}
	0 & 0\\[2mm]
	c & ic
	\end{pmatrix}.\]

	For the algebra  $E_{2}$: $e_1e_1=e_1, \   e_2e_2=e_1$ to find the matrix form of the Rota-Baxter operator of weight $1$ we should solve the following system of equations:
	
	\[\left\{\begin{array}{llllll}
	b^2=a^2+a, \\[1mm]
	(2a+1)b=0, \\[1mm]
	c^2+d^2=(2d+1)a, \\[1mm]
	(2d+1)b=0, \\[1mm]
	bd=ab, \\[1mm]
	(b+c)b=0.
	\end{array}
	\right.\]
	
	This system of equations has the following solutions:
	
	\[\begin{pmatrix}
	0 & 0\\[2mm]
	c & i c
	\end{pmatrix}, \quad
	\begin{pmatrix}
	0 & 0\\[2mm]
	c & -i c
	\end{pmatrix}, \quad
	\begin{pmatrix}
	-\frac{1}{2} & \frac{i}{2}\\[2mm]
	-\frac{i}{2} & -\frac{1}{2}
	\end{pmatrix}, \quad
	\begin{pmatrix}
	-\frac{1}{2} & -\frac{i}{2}\\[2mm]
	\frac{i}{2} & -\frac{1}{2}
	\end{pmatrix},\]
	
	\[\begin{pmatrix}
	-1 & 0\\[2mm]
	c & -1+ic
	\end{pmatrix}, \quad
	\begin{pmatrix}
	-1 & 0\\[2mm]
	c & -1-ic
	\end{pmatrix}.\]

	\item For the algebra $E_{3}: e_1e_1=e_1+e_2, \ \  e_2e_2=-e_1-e_2$ to find the matrix form of the Rota-Baxter operator of weight $0$ we should solve the next system of equations:
	
	\[\left\{\begin{array}{llllll}
	a^2-b^2=2a(a+c), \\[1mm]
	a^2-b^2=2a(b+d), \\[1mm]
	c^2-d^2=-2d(a+c), \\[1mm]
	c^2-d^2=-2d(b+d), \\[1mm]
	ac-bd=(c-b)(a+c), \\[1mm]
	ac-bd=(c-b)(b+d).
	\end{array}
	\right.\]
	
	This system of equations has the following solutions:
	
	\[\begin{pmatrix}
	a & a\\[2mm]
	-a & -a
	\end{pmatrix}, \qquad
\begin{pmatrix}
	a & -a\\[2mm]
	-a & a
	\end{pmatrix}.\]

	For the algebra $E_{3}: e_1e_1=e_1+e_2, \  e_2e_2=-e_1-e_2$ to find the matrix form of the Rota-Baxter operator of weight $1$ we should solve the following system of equations:
	
	\[\left\{\begin{array}{llllll}
	a^2-b^2=(2a+1)(a+c), \\[1mm]
	a^2-b^2=(2a+1)(b+d), \\[1mm]
	c^2-d^2=-(2d+1)(a+c), \\[1mm]
	c^2-d^2=-(2d+1)(b+d), \\[1mm]
	ac-bd=(c-b)(a+c), \\[1mm]
	ac-bd=(c-b)(b+d).
	\end{array}
	\right.\]
	
	This system of equations has the following solutions:
\[ \begin{pmatrix}
	-1+b & b\\
	-b & -1-b
	\end{pmatrix}, \qquad \begin{pmatrix}
	-1-b & b\\
	b & -1-b
	\end{pmatrix}, \qquad  
	 \begin{pmatrix}
	b & b\\
	-b & -b
	\end{pmatrix}, \qquad  
	\begin{pmatrix}
	-b & b\\
	b & -b
	\end{pmatrix}.\]

	\item For the algebra $E_4: e_1e_1=e_2$ to find the matrix form of the Rota-Baxter operator of weight $0$ we should solve the following system of equations:
	
	\[\left\{\begin{array}{llllll}
	a^2=2ad, \\[1mm]
	2ac=0, \\[1mm]
	ac=cd, \\[1mm]
	c^2=0.
	\end{array}
	\right.\]
	
	This system of equations has the next solutions:

	\[\begin{pmatrix}
	0 & b\\
	0 & d
	\end{pmatrix}, \quad
	\begin{pmatrix}
	a & b\\
	0 & \frac{a}{2}
	\end{pmatrix}.\]

	For the algebra $E_4: e_1e_1=e_2$ to find the matrix form of the Rota-Baxter operator of weight $1$ we should solve the following system of equations:

	\[\left\{\begin{array}{llllll}
	a^2=(2a+1)d, \\[1mm]
	(2a+1)c=0, \\[1mm]
	ac=cd, \\[1mm]
	c^2=0.
	\end{array}
	\right.\]
	
	This system of equations has the following solution:
	
	\[\begin{pmatrix}
	a & b \\
	0 & \frac{a^2}{1+2a}
	\end{pmatrix},\]
	where $a\neq -\frac{1}{2}$. The system has not solution when $a=-\frac{1}{2}$.

	\item For the algebra $E_5(x,y): e_1e_1=e_1+x e_2, \ \ e_2e_{2}=y e_1+e_2, \ \ 1-x y \neq 0, \ x, y\in {\mathbb{C}}$, to find  the matrix form of the Rota-Baxter operator of weight $1$ we should solve the following system of equations:
	
	\begin{equation}\label{rb5}
	\left\{\begin{array}{llllll}
	a^2+b^2 y=(2a+1)(a+xc), \\[1mm]
	a^2x+b^2=(2a+1)(b+xd), \\[1mm]
	c^2+d^2 y=(2d+1)(ay+c), \\[1mm]
	c^2x+d^2=(2d+1)(by+d), \\[1mm]
	ac+bdy=c(a+cx)+b(ay+c), \\[1mm]
	acx+bd=c(b+dx)+b(by+d).
	\end{array}
	\right.
	\end{equation}
	
	To find the solution of  system \eqref{rb5} consider the following:
	
	\textbf{Case 1.} Let $a=0$. Then from \eqref{rb5}

	\begin{equation}\label{rb51}
	\left\{\begin{array}{llllll}
	b^2 y=xc, \\[1mm]
	b^2=b+xd, \\[1mm]
	c^2+d^2 y=(2d+1)c, \\[1mm]
	c^2x+d^2=(2d+1)(by+d), \\[1mm]
	bdy=c^2x+bc, \\[1mm]
	bd=c(b+dx)+b(by+d).
	\end{array}
	\right.
	\end{equation}

	From this system consider:
	
	\textbf{Case 1.1.} Let $b=0$. Then
	
	\[\left\{\begin{array}{llllll}
	0=xc, \\[1mm]
	0=xd, \\[1mm]
	c^2+d^2 y=(2d+1)c, \\[1mm]
	c^2x+d^2=(2d+1)d, \\[1mm]
	0=c^2x, \\[1mm]
	0=cdx.
	\end{array}
	\right.\]
	
	In the case when $x=0$ we have:
	
	\[ \begin{pmatrix}
	0 & 0\\
	1 & 0
	\end{pmatrix}, \quad
	\begin{pmatrix}
	0 & 0\\
	c_1 & -1
	\end{pmatrix}, \quad
	\begin{pmatrix}
	0 & 0 \\
	c_2 & -1
	\end{pmatrix},\]
	where $c_{1,2}=\frac{-1 \pm \sqrt{1-4y}}{2}$.
	
	When $x\neq 0$ we have the trivial solution of the system $a=b=c=d=0$.
	
	\textbf{Case 1.2.} $b \neq 0$. Then in \eqref{rb51} consider:
	
	\textbf{Case 1.2.1.} $c=0$. Then we get:
	
	\[\left\{\begin{array}{llllll}
	b^2 y=0, \\[1mm]
	b^2=b+xd, \\[1mm]
	d^2 y=0, \\[1mm]
	d^2=(2d+1)(by+d), \\[1mm]
	bdy=0, \\[1mm]
	b^2y=0.
	\end{array}
	\right.\]
	
	Since $b \neq 0$, we have $y=0$. Then we have:
	
	\[\begin{pmatrix}
	0 & 1\\
	0 & 0
	\end{pmatrix}, \quad
	\begin{pmatrix}
	0 & b_1\\
	0 & -1
	\end{pmatrix}, \quad
	\begin{pmatrix}
	0 & b_2\\
	0 & -1
	\end{pmatrix},\]
	where $b_{1,2}=\frac{1 \pm \sqrt{1-4x}}{2}$.

	\textbf{Case 1.2.2.} Let $c \neq 0$. Then from \eqref{rb51} consider:
	
	\textbf{Case 1.2.2.1.} When $x=0$, then $y=0$. And we get the system
	
	\[\left\{\begin{array}{llllll}
	b^2=b, \\[1mm]
	c^2=(2d+1)c, \\[1mm]
	d^2+d=0, \\[1mm]
	bc=0,
	\end{array}
	\right.\]
	which has not solution.
	
	\textbf{Case 1.2.2.2.} Let $x \neq 0$. Since $b \neq 0, \ c \neq 0$, then $y \neq 0$. Then from \eqref{rb51} consider:

	\textbf{Case 1.2.2.2.1.} $d=0$. And so we have the system
	
	\[\left\{\begin{array}{llllll}
	b^2 y=c x, \\[1mm]
	b^2=b, \\[1mm]
	c^2=c, \\[1mm]
	c^2 x=by, \\[1mm]
	c^2 x+bc=0, \\[1mm]
	b^2 y+bc=0,
	\end{array}
	\right.\]
	which has solution $b=c=1, x=y=-1$. But this solution does not satisfy the condition $xy \neq 1$.
	
	\textbf{Case 1.2.2.2.2.} $d \neq 0$. Then from \eqref{rb51} consider:
	
	\textbf{Case 1.2.2.2.2.1.} When $d=-\frac{1}{2}$ we get a solution, which  also contradicts the condition $xy \neq 1$.
	
	\textbf{Case 1.2.2.2.2.2.} When $d \neq -\frac{1}{2}$. One have to solve the system:

	\[	\left\{\begin{array}{llllll}
	b^2 y=xc, \\[1mm]
	b^2=b+xd, \\[1mm]
	c^2+d^2 y=2cd+c, \\[1mm]
	c^2x=d^2+2bdy+by+d, \\[1mm]
	bdy=c^2x+bc, \\[1mm]
	bc+cdx+b^2y=0.
	\end{array}
	\right.\]
	
The system has the next solution: $b=-1, \ c=1, \ d=-2, \ x=-1, \ y=-1$. But this solution also contradicts the condition $xy \neq 1$.
		
	\textbf{Case 2.} Let $a \neq 0$. Then from \eqref{rb5}:
	
	\textbf{Case 2.1.} For $a=-\frac{1}{2}$ we get:
	
	\[\left\{\begin{array}{llllll}
	\frac{1}{4}+b^2 y=0, \\[1mm]
	\frac{x}{4}+b^2=0, \\[1mm]
	c^2+d^2 y=(2d+1)(-\frac{y}{2}+c), \\[1mm]
	c^x +d^2=(2d+1)(by+d), \\[1mm]
	-\frac{c}{2}+bdy=c(-\frac{1}{2}+cx)+b(-\frac{y}{2}+c), \\[1mm]
	-\frac{cx}{2}+bd=c(b+xd)+b(by+d).
	\end{array}
	\right.\]
	
	From the first two equations of the system, we can see that $xy=1$, which contradicts the condition $xy \neq 1$.
	
	\textbf{Case 2.2.} Let $a \neq -\frac{1}{2}$. Then from \eqref{rb5}:
	
	\textbf{Case 2.2.1.} When $b=0$, we get
	
	\[\left\{\begin{array}{llllll}
	a^2+2acx+a+cx=0, \\[1mm]
	a^2x-2axd-xd=0, \\[1mm]
	c^2+d^2 y=2ady+2cd+ay+c, \\[1mm]
	c^2x=d^2+d, \\[1mm]
	c^2x=0, \\[1mm]
	acx-cdx=0.
	\end{array}
	\right.\]
	
	For this system, when $c=0$, we have the following solutions:
	
	$\begin{pmatrix*}[r]
	-1 & 0\\
	0 & 0
	\end{pmatrix*}$,
	which corresponds to the evolution algebra $E(0,0)$.
	
	$\begin{pmatrix*}[r]
	-1 & 0\\
	0 & -1
	\end{pmatrix*}$,
	which corresponds to the evolution algebra $E(x,y)$.
	
	When $c \neq 0$, it implies that $x=0$. Thus we have the next solutions that correspond to $E(0,y)$:
	
	\[\begin{pmatrix}
	-1 & 0 \\
	-1 & -1
	\end{pmatrix}, \quad
	\begin{pmatrix}
	-1 & 0\\
	c_1 & 0
		\end{pmatrix}, \quad
	\begin{pmatrix}
	-1 & 0\\
	c_2 & 0
	\end{pmatrix},\]
	where $c_{1,2}=\frac{1 \pm \sqrt{1-4y}}{2}$ and $y\ne 0$.

	\textbf{Case 2.2.2.} Let $b \neq0$. Then from \eqref{rb5} consider:
	
	\textbf{Case 2.2.2.1.} $c=0$. Then we have
	
	\[\left\{\begin{array}{llllll}
	b^2 y=a^2+a, \\[1mm]
	a^2x+b^2=2ab+2adx+b+xd, \\[1mm]
	d^2 y=2ady+ay, \\[1mm]
	d^2+2bdy+by+d=0, \\[1mm]
	bdy=aby, \\[1mm]
	b^2y=0.
	\end{array}
	\right.\]
	
	From the last equation of the system, it implies that $y=0$. Thus we have the following solutions corresponding to
	the evolution algebra $E(x,0)$:
	
	\[\begin{pmatrix}
	-1 & b_{1,2}\\[2mm]
	0 & 0
	\end{pmatrix}, \quad
	\begin{pmatrix}
	-1 & -1\\[2mm]
	0 & -1
	\end{pmatrix},\]
	where $b_{1,2}=\frac{-1 \pm \sqrt{1-4x}}{2}$.

	\textbf{Case 2.2.2.2.} Let $c\neq 0$. Then from \eqref{rb5} consider:
	
	\textbf{Case 2.2.2.2.1.} $d= 0$. Then we have
	
	\[\left\{\begin{array}{llllll}
	b^2 y=a^2+2acx+a+cx, \\[1mm]
	a^2x+b^2=2ab+b, \\[1mm]
	c^2=ay+c, \\[1mm]
	c^2x=by, \\[1mm]
	c^2x+aby+bc=0, \\[1mm]
	acx=bc+b^2y.
	\end{array}
	\right.\]

	When $x=0$, it is easy to see that the system of equations has not solution.
	
	From $x\neq 0$ implies that $y\neq 0$. Then

	\textbf{Case 2.2.2.2.2.} Let $d\neq 0$. When $d=-\frac{1}{2}$ we get a solution, which  also contradicts the condition $xy \neq 1$. Then from \eqref{rb5} consider:
	
	\textbf{Case 2.2.2.2.2.1.} $x= 0$. Then we have

	\[\left\{\begin{array}{llllll}
	b^2 y=a^2+a, \\[1mm]
	b^2=2ab+b, \\[1mm]
	c^2+d^2 y=2ady+2cd+ay+c, \\[1mm]
	d^2+2bdy+by+d=0, \\[1mm]
	bdy=aby+bc, \\[1mm]
	bc+b^2y=0.
	\end{array}
	\right.\]
	
	When $y=0$ and $y=-\frac{1}{4}$ one can check that the system has not solutions.
	
	In the case, when $y\neq 0$ and $y\neq \frac{1}{4}$, we have the next solutions of the system:
	
	\[	\begin{pmatrix}
	\frac{1-4y+\sqrt{1-4y}}{8y-2} & -\frac{1}{\sqrt{1-4y}}\\[2mm]
	\frac{y}{\sqrt{1-4y}} & \frac{1-4y-\sqrt{1-4y}}{8y-2}
		\end{pmatrix}, \quad
	\begin{pmatrix}
	\frac{1-4y-\sqrt{1-4y}}{8y-2} & \frac{1}{\sqrt{1-4y}}\\[2mm]
	-\frac{y}{\sqrt{1-4y}} & \frac{1-4y+\sqrt{1-4y}}{8y-2}
		\end{pmatrix}.\]

	\textbf{Case 2.2.2.2.2.2.} When $x\neq 0$. Then from \eqref{rb5} consider:
	
	\textbf{Case 2.2.2.2.2.2.1.} $y=0$. Then we get
	
	\[\left\{\begin{array}{llllll}
	a^2+2acx+a+cx=0, \\[1mm]
	a^2 x+b^2=2ab+2adx+b+dx, \\[1mm]
	c^2=2cd+c, \\[1mm]
	c^2x=d^2+d, \\[1mm]
	c^2x+bc=0, \\[1mm]
	acx=bc+cdx.
	\end{array}
	\right.\]
	
	The system of equations has the following solutions:
	
	\[\begin{pmatrix}
	\frac{1-4x+\sqrt{1-4x}}{8x-2} & -\frac{x}{\sqrt{1-4x}}\\[2mm]
	\frac{1}{\sqrt{1-4x}} & \frac{1-4x-\sqrt{1-4x}}{8x-2}
	\end{pmatrix}, \quad
\begin{pmatrix}
	\frac{1-4x-\sqrt{1-4x}}{8x-2} & \frac{x}{\sqrt{1-4x}}\\[2mm]
	-\frac{1}{\sqrt{1-4x}} & \frac{1-4x+\sqrt{1-4x}}{8x-2}
	\end{pmatrix},\]
	where $x\neq \frac{1}{4}$.
		
	\textbf{Case 2.2.2.2.2.2.2.}  When $y\neq 0$, then we should solve \eqref{rb5} for all non-zero unknowns.
	
	\begin{equation}\label{eqneq}
	\left\{\begin{array}{llllll}
	a^2+b^2 y=(2a+1)(a+xc), \\[1mm]
	a^2x+b^2=(2a+1)(b+xd), \\[1mm]
	c^2+d^2 y=(2d+1)(ay+c), \\[1mm]
	c^2x+d^2=(2d+1)(by+d), \\[1mm]
	ac+bdy=c(a+cx)+b(ay+c), \\[1mm]
	acx+bd=c(b+dx)+b(by+d).
	\end{array}
	\right.
	\end{equation}

	From the second and the third equations of the system $(a^2-d-2ad)x=b+2ab-b^2, \ (d^2-a-2ad)y=c-c^2+2cd$, consider the next cases:
	
	\underline{\textbf{Case A.}} Let $a^2-d-2ad=0, b+2ab-b^2=0, d^2-a-2ad=0,c-c^2+2cd=0$. From which we will get the next solutions:

	$a) \ \ a=b=c=d=-1$;

	$b) \ \ a=\frac{-3+i \sqrt{3}}{6}, b=\frac{i}{\sqrt{3}}, c=-\frac{i}{\sqrt{3}}, d=\frac{-3-i \sqrt{3}}{6}$;

	$c) \ \ a=\frac{-3-i \sqrt{3}}{6}, b=-\frac{i}{\sqrt{3}}, c=\frac{i}{\sqrt{3}}, d=\frac{-3+i \sqrt{3}}{6}$.
	
	$a=b=c=d=-1$ will be solution of the system of equations \eqref{eqneq} if and only if  $x=y=1$, which contradicts the condition $xy\neq 1$.
	
	$a=\frac{-3+i \sqrt{3}}{6}, \ \ b=\frac{i}{\sqrt{3}}, \ \ c=-\frac{i}{\sqrt{3}}, d=\frac{-3-i \sqrt{3}}{6}$ and 	
	$a=\frac{-3-i \sqrt{3}}{6}, \ \ b=-\frac{i}{\sqrt{3}}, \ \ c=\frac{i}{\sqrt{3}}, \ \ d=\frac{-3+i \sqrt{3}}{6}$ will be solutions of the system of equations \eqref{eqneq} if and only if $x=1-y$. From the condition $xy\neq 1$ we get $x\neq \frac{1\pm i\sqrt{3}}{2}$.
	
\underline{\textbf{Case B.}} Let $a^2-d-2ad=0, \  b+2ab-b^2=0, \  y=\frac{c-c^2+2cd}{d^2-a-2ad}$. Then substituting
$d=\frac{a^2}{2a+1}, \ b=2a+1, \ y=\frac{c-c^2+2cd}{d^2-a-2ad}$, we get the next:

\[\left\{\begin{array}{lll}
\begin{array}{cc}
a+a^2+\frac{(1+2a)^3(1+2a^2-2a(-1+c)-c)c}{a(1+4a+6a^2+3a^3)}+cx+2acx=0;
\end{array}\\[4mm]
\begin{array}{ccc}
\frac{(1+2a)^2\left((1+2a)(1+\frac{2a^2}{1+2a})c(-1-\frac{2a^2}{1+2a}+c)+
	\frac{a(1+4a+6a^2+3a^3)(\frac{a^4}{(1+2a)^2})+\frac{a^2}{1+2a}-c^2x}{(1+2a)^2}\right)}{a(1+4a+6a^2+3a^3)}=0;
\end{array}\\[6mm]
\begin{array}{ccc}
\frac{c(2a^3+c(1+x)+a(1+c(4+3x))+a^2(3+c(4+3x)))}{1+3a+3a^2}=0;
\end{array}\\ [4mm]
\begin{array}{ccc}
c(-1-2a+\frac{(1+2a)^3(1+2a^2-2a(-1+c)-c)}{a(1+4a+6a^2+3a^3)}+\frac{a(1+a)x}{1+2a})=0,
\end{array}
\end{array}\right.\]
where $a\neq-\frac{1}{2}, \ a\neq -1, \ a\neq \frac{-3+i \sqrt{3}}{6}, a\neq \frac{-3-i \sqrt{3}}{6}$.

From which we will get the following solution:

 $x_1=-\frac{(1+2a)^2}{a^2},  \quad c_1=\frac{a^3}{(1+2a)^2}$, \quad
 $x_2=-\frac{(1+2a)^2}{(1+a)^2}, \quad \ c_2=\frac{1+3a+3a^2+a^3}{(1+2a)^2}$.

Thus we will get the solution of the system in the next form:

$b=1+2a$,
$c_1=\frac{a^3}{(1+2a)^2}$,
$d=\frac{a^2}{1+2a}$,
$x_1=-\frac{(1+2a)^2}{a^2}$,
$y_1=-\frac{a^2}{(1+2a)^2}$,

$b=1+2a$,
$c_2=\frac{1+3a+3a^2+a^3}{(1+2a)^2}$,
$d=\frac{a^2}{1+2a}$,
$x_2=-\frac{(1+2a)^2}{(1+a)^2}$,
$y_2=-\frac{(1+a)^2}{(1+2a)^2}$.

But these solutions contradict the condition $xy\neq 1$.

\underline{\textbf{Case C.}} Let $x=\frac{b+2ab-b^2}{a^2-d-2ad}, \ c-c^2+2cd=0, \ d^2-a-2ad=0$.
Similar to the Case B we get the next solutions of the system:

$a=\frac{d^2}{1+2d}$,
$b_1=\frac{d^3}{(1+2d)^2}$,
$c=1+2d$,
$x_1=-\frac{d^2}{(1+2d)^2}$,
$y_1=-\frac{(1+2d)^2}{d^2}$,

$a=\frac{d^2}{1+2d}$,
$b_2=\frac{1+3d+3d^2+d^3}{(1+2d)^2}$,
$c=1+2d$,
$x_2=-\frac{(1+d)^2}{(1+2d)^2}$,
$y_2=-\frac{(1+2d)^2}{(1+d)^2}$.

Also these solutions contradict the condition $xy\neq 1$.

Thus we cannot take these solutions as a solution of  system \eqref{eqneq}.

\underline{\textbf{Case D.}} Let $x=\frac{b+2ab-b^2}{a^2-d-2ad}, \ y=\frac{c-c^2+2cd}{d^2-a-2ad}$.

In this case we have the next solution of  system \eqref{eqneq}:

$a=-1-d$,
$b=-\frac{d(1+d)}{c}$,
$x=\frac{d(1+d)(c+2cd-d(1+d))}{c^2(1+3d+3d^2)}$,
$y=\frac{c(1-c+2d)}{c^2(1+3d+3d^2)}$,

\noindent where $d\neq -1, \ d\neq -\frac{3\pm i \sqrt{3}}{6}, \ c\neq \frac{d(1+d)}{1+2d}, \ c\neq 1+2d$.

	For the algebra $E_5(x,y): e_1e_1=e_1+x e_2, \ \ e_2e_{2}=y e_1+e_2, \ \ 1-x y \neq 0, \ x, y \in {\mathbb{C}}$,
to find the matrix form of the Rota-Baxter operator of weight $0$ we should solve the following system of equations:
	
	\begin{equation}\label{rb50}
	\left\{\begin{array}{llllll}
	b^2 y=a^2+2acx, \\[1mm]
	a^2x+b^2=2ab+2adx, \\[1mm]
	c^2+d^2 y=2ady+2cd, \\[1mm]
	c^2x=2bdy+d^2, \\[1mm]
	bdy=c^2x+aby+bc, \\[1mm]
	acx=bc+cdx+b^2y.
	\end{array}
	\right.
	\end{equation}
	
	One can check that this system of equations has only the trivial solution, in the case when one of unknowns of $a, b, c, d$ will equal to $0$.
	
	Also, when $x=y=0$, we have only the trivial solution of the system.
	
	For the algebra $E_5(0,\frac{1}{4})$,   system  \eqref{rb50} has the following solution:
	
	\[\begin{pmatrix}
	a & 2a\\[2mm]
	-\frac{a}{2} & -a
	\end{pmatrix}.\]

	And for the algebra $E_5(\frac{1}{4},0)$,  system  \eqref{rb50} has the following solution:
	
	\[\begin{pmatrix}
	a & \frac{a}{2}\\[2mm]
	-2a & -a
	\end{pmatrix}.\]
	
	Consider for system \eqref{rb50}, that $a\neq 0, b\neq 0, c\neq 0, d\neq 0$, and for the algebra $E(x,y)$, let $x\neq 0, y \neq0$.
	
	From the first and the fifth equations of  system \eqref{rb50}, let $b^2 y-a^2\neq o, \ \ c^2x+aby+bc \neq 0$ (otherwise we get a contradiction). Then,
	
	\[c=\frac{b^2 y-a^2}{2ax}, \quad  d=\frac{c^2x+aby+bc}{by}.\]
	
	After replacing these values in the second equation of system \eqref{rb50}, we have
	
	\[x=\frac{-a^4+2a^3b-2a^2b^2y-b^4y^2}{2a^3by}.\]
	
	Then substituting these values in the $3,4,6$-th equations of system \eqref{rb50} we get
	
	\[y=\frac{-a^2+2ab}{3b^2}.\]
	
	Then substituting the value of $y$ in the value of $x$, we get
	
	\[x=\frac{(2a-b)b}{3a^2}.\]
	
	For the parameters $x, y$ we should check that $xy \neq 1$. Substituting the values of $x$ and $y$, we have
	
	\[xy=-\frac{(a-2b)(2a-b)}{9ab}\neq 1.\]
	
	Then, we get $a \neq -b$.
	
	Substituting the values of $x, y$ in the values of $c,d$, we get $d=-a, \ c=-\frac{a^2}{b}$.
	
	Thus, we have the following  solution of system \eqref{rb50}:
	
	\[\begin{pmatrix}
	a & b \\[2mm]
	-\frac{a^2}{b} & -a
	\end{pmatrix},
	\]
	where $a \neq -b$.
	
	\item For the algebra $E_6(x): e_1e_1=e_2, \ e_2e_{2}=e_1+ x e_2, \ x \in {\mathbb{C}}$, to find the matrix form of the Rota-Baxter operator of weight $1$ we should solve the following system of equations:
	
	\begin{equation}\label{rb61}
	\left\{\begin{array}{llllll}
	b^2=(2a+1)c, \\[1mm]
	a^2+b^2 x=(2a+1)d, \\[1mm]
	d^2=(2d+1)(a+cx), \\[1mm]
	c^2+d^2 x=(2d+1)(b+dx), \\[1mm]
	bd=ab+c^2+bcx, \\[1mm]
	ac=cd+b^2.
	\end{array}
	\right.
	\end{equation}
	
	\textbf{Case 1.} Let $a=0$. Then
	\[
	\left\{\begin{array}{llllll}
	b^2=c, \\[1mm]
	b^2 x=d, \\[1mm]
	d^2=2cdx+cx, \\[1mm]
	c^2+d^2 x=2bd+2d^2x+b+dx, \\[1mm]
	bd=c^2+bcx, \\[1mm]
	cd+b^2=0.
	\end{array}
	\right.
	\]
	
	\textbf{Case 1.1.} Consider $b=0$. Then it is easy to see that the system has only the trivial solution.
	
	\textbf{Case 1.2.} In the case when $b \neq 0$, one can check that the system has not solutions.
	
	\textbf{Case 2.} Let $a\neq 0$ and $a \neq -\frac{1}{2}$ (in the case when $a=-\frac{1}{2}$, it is easy to see that \eqref{rb61} has not solution).
	
	\textbf{Case 2.1.} Let $b=0$. Then from \eqref{rb61},
	
	\[
	\left\{\begin{array}{llllll}
	(2a+1)c=0, \\[1mm]
	a^2=(2a+1)d, \\[1mm]
	d^2=(2d+1)(a+cx), \\[1mm]
	c^2+d^2 x=2d^2x+dx, \\[1mm]
	c^2=0, \\[1mm]
	ac=cd.
	\end{array}
	\right.
	\]
	
	And this system has the following solutions:
	
	For the algebra $E_6(x)$
	
	\[\begin{pmatrix}
	-1 & 0\\[2mm]
	0 & -1
	\end{pmatrix},
	\]
	
	and for the algebra $E_6(0)$
	
	\[\begin{pmatrix}
	\frac{-3+i\sqrt{3}}{6} & 0\\[2mm]
	0 & \frac{-3-i\sqrt{3}}{6}
	\end{pmatrix},\quad
	\begin{pmatrix}
	\frac{-3-i\sqrt{3}}{6} & 0\\[2mm]
	0 & \frac{-3+i\sqrt{3}}{6}
	\end{pmatrix}.
	\]

	\textbf{Case 2.2.} Let $b \neq 0$.
	
	\textbf{Case 2.2.1.} Let $c=0$. Then we get $b=0$, which contradicts the \textbf{Case 2.2.}. Thus we do not have solution.
	
	\textbf{Case 2.2.2.} Let $c \neq 0$. Then for the case, when $d=0$, the system has not solution. Consider $d\neq 0$.
	
	\textbf{Case 2.2.2.1.} Consider $x=0$. Then from \eqref{rb61} we get the following system of equations:
	
	\[
	\left\{\begin{array}{llllll}
	b^2=(2a+1)c, \\[1mm]
	a^2=(2a+1)d, \\[1mm]
	d^2=2ad+a, \\[1mm]
	c^2=2bd+b, \\[1mm]
	bd=ab+c^2, \\[1mm]
	ac=cd+b^2.
	\end{array}
	\right.
	\]
	
	This system has the following solutions:
	
	\[\begin{pmatrix}
	\frac{-3+i\sqrt{3}}{6} & -\frac{i}{\sqrt{3}}\\[2mm]
	\frac{i}{\sqrt{3}} & \frac{-3-i\sqrt{3}}{6}
	\end{pmatrix},\quad
	\begin{pmatrix}
	\frac{-3-i\sqrt{3}}{6} & \frac{i}{\sqrt{3}}\\[2mm]
	-\frac{i}{\sqrt{3}} & \frac{-3+i\sqrt{3}}{6}
	\end{pmatrix},\quad
	\begin{pmatrix}
	\frac{-3-i\sqrt{3}}{6} & -\frac{\sqrt[6]{-1}}{\sqrt{3}}\\[2mm]
	\frac{(\sqrt[6]{-1})^5}{\sqrt{3}} & \frac{-3+i\sqrt{3}}{6}
	\end{pmatrix},
	\]

	\[\begin{pmatrix}
	\frac{-3+i\sqrt{3}}{6} & \frac{\sqrt[6]{-1}}{\sqrt{3}}\\[2mm]
	-\frac{(\sqrt[6]{-1})^5}{\sqrt{3}} & \frac{-3-i\sqrt{3}}{6}
	\end{pmatrix},\quad
	\begin{pmatrix}
	\frac{-3+i\sqrt{3}}{6} & \frac{(\sqrt[6]{-1})^5}{\sqrt{3}}\\[2mm]
	\frac{\sqrt[6]{-1}}{\sqrt{3}} & \frac{-3-i\sqrt{3}}{6}
	\end{pmatrix},\quad
	\begin{pmatrix}
	\frac{-3-i\sqrt{3}}{6} & -\frac{(\sqrt[6]{-1})^5}{\sqrt{3}}\\[2mm]
	\frac{\sqrt[6]{-1}}{\sqrt{3}} & \frac{-3+i\sqrt{3}}{6}
	\end{pmatrix}.
	\]

	\textbf{Case 2.2.2.2.} Consider $x \neq 0$. Then we get the following solution of  system \eqref{rb61} for the algebra $E_6(\frac{-b^3-c^3}{bc^2})$:

	\[\begin{pmatrix}
	\frac{b^2-c}{2c} & b\\[2mm]
	c & \frac{-b^2-c}{2c}
	\end{pmatrix},
	\]
	where the parameters $b,c$ are solutions of the following system of equations

\[\begin{cases}
	\frac{b^6+5b^3c^3+4c^6}{c}&=\frac{c(b^3+c^3)}{b}, \\[1mm]
	\quad \frac{b^4}{c}+4bc^2&=c.
	\end{cases}\]

	For the algebra $E_6(x): e_1e_1=e_2, \ e_2e_{2}=e_1+ x e_2, \ x \in {\mathbb{C}}$, to find the matrix form of the Rota-Baxter operator of weight $0$ we should solve the following system of equations:
	
	\begin{equation*}%\label{rb60}
	\left\{\begin{array}{llllll}
	b^2=2ac, \\[1mm]
	a^2+b^2 x=2ad, \\[1mm]
	d^2=2ad+2cdx, \\[1mm]
	c^2=2bd+d^2x, \\[1mm]
	bd=ab+c^2+bcx, \\[1mm]
	ac=cd+b^2.
	\end{array}
	\right.
	\end{equation*}
	
	Which will have the following solution for the algebra $E_6(-\frac{3b^2}{4c^2})$:
	
		\[\begin{pmatrix}
	\frac{b^2}{2c} & b\\[2mm]
	c & \frac{-b^2}{2c}
	\end{pmatrix},
	\]
	where the parameters $b,c$ are solutions of the following system of equations

\[	\begin{cases}
	\frac{3b^6}{c}+16b^3c^2+16c^5&=0, \\[1mm]
	\qquad \qquad \ \ \frac{b^3}{c}+4c^2&=0.
	\end{cases}\]

\end{itemize}

Thus we have proved the following theorem, which gives all matrices form of the Rota-Baxter operators on 2-dimensional complex evolution algebras.

\begin{thm}
	The matrices of the Rota-Baxter operators on the two-dimensional complex evolution algebras are given in the next table, with parameters $a,b,c,d,x,y \in {\mathbb{C}}$.

	\

\begin{tabular}{|c|c|}
		\hline
		Evolution Algebra &  Matrices of RBOs of weight 0 \\
		 {} & on the evolution algebra \\
		\hline \hline
		$E_1$	& $\begin{pmatrix}
		0 & b \\[2mm]
		0 & d
		\end{pmatrix}$ \\
		\hline
		$E_2$	& $\begin{pmatrix}
		0 & 0\\[2mm]
		c & -ic
		\end{pmatrix}, \quad
		\begin{pmatrix}
		0 & 0\\[2mm]
		c & ic
		\end{pmatrix}$ \\
		\hline
		$E_3$	& $\begin{pmatrix}
		a & a\\[2mm]
		-a & -a
		\end{pmatrix}, \quad
		\begin{pmatrix}
		a & -a\\[2mm]
		-a & a
		\end{pmatrix}$ \\
		\hline
		$E_4$	& $\begin{pmatrix}
		0 & b\\[2mm]
		0 & d
		\end{pmatrix}, \quad
		\begin{pmatrix}
		a & b\\[2mm]
		0 & \frac{a}{2}
		\end{pmatrix}$ \\
		\hline
		$E_5(\frac{1}{4},0)$	& 	$
		\begin{pmatrix}
		a & \frac{a}{2}\\[2mm]
		-2a & -a
		\end{pmatrix}$ \\
		\hline
		$E_5(0,\frac{1}{4})$	&
		$\begin{pmatrix}
		a & 2a\\[2mm]
		-\frac{a}{2} & -a
		\end{pmatrix}$  \\
		\hline
		$E_5\big(\frac{(2a-b)b}{3a^2},\frac{-a^2+2ab}{3b^2}\big)$ & {}\\
		$a\neq2b, b\neq 2a$, $a \neq -b$ 	&
		$\begin{pmatrix}
		a & b \\[2mm]
		-\frac{a^2}{b} & -a
		\end{pmatrix}$  \\
    	$a\neq 0, b\neq 0$	& {}\\
		\hline
		$E_6(-\frac{3b^2}{4c^2})$	&
		$\begin{pmatrix}
		\frac{b^2}{2c} & b\\[2mm]
		c & \frac{-b^2}{2c}
		\end{pmatrix}$ \\
		$b\neq 0, \ c\neq 0$& where the parameters $b,c$ are solutions\\
		&  of the following system of equations\\	
		& $	\begin{cases}
		\frac{3b^6}{c}+16b^3c^2+16c^5&=0, \\[1mm]
		\qquad \qquad \ \ \frac{b^3}{c}+4c^2&=0.
		\end{cases}$\\
		\hline
	\end{tabular}

\newpage
	
\noindent \begin{longtable}{|c|c|}
		\hline
		Evolution Algebra & Matrices of RBOs of weight 1 \\
		{} &  on the evolution algebra \\
		\hline \hline
		$E_1$	& $\begin{pmatrix}
		-1 & 0\\[2mm]
		0 & d
		\end{pmatrix}, \quad
		\begin{pmatrix}
		0 & 0\\[2mm]
		0 & d
		\end{pmatrix}$\\
		\hline
		$E_2$	&  $\begin{pmatrix}
		0 & 0\\[2mm]
		c & i c
		\end{pmatrix}, \quad
		\begin{pmatrix}
		0 & 0\\[2mm]
		c & -i c
		\end{pmatrix}$,\\
		& $\begin{pmatrix}
		-\frac{1}{2} & \frac{i}{2}\\[2mm]
		-\frac{i}{2} & -\frac{1}{2}
		\end{pmatrix}, \quad
		\begin{pmatrix}
		-\frac{1}{2} & -\frac{i}{2}\\[2mm]
		\frac{i}{2} & -\frac{1}{2}
		\end{pmatrix}$,
		\\
		& $\begin{pmatrix}
		-1 & 0\\[2mm]
		c & -1+ic
		\end{pmatrix}, \quad
		\begin{pmatrix}
		-1 & 0\\[2mm]
		c & -1-ic
		\end{pmatrix}$ \\
		    {} & {} \\
		\hline
		$E_3$	& $\begin{pmatrix}
		-1+b & b\\[2mm]
		-b & -1-b
		\end{pmatrix}, \
		\begin{pmatrix}
		-1-b & b\\[2mm]
		b & -1-b
		\end{pmatrix}$, \\ [6mm]
		& $\begin{pmatrix}
		b & b\\[2mm]
		-b & -b
		\end{pmatrix}, \quad
		\begin{pmatrix}
		-b & b\\[2mm]
		b & -b
		\end{pmatrix}$ \\
 {} & {} \\
		\hline
		$E_4$	& $\begin{pmatrix}
		a & b \\[2mm]
		0 & \frac{a^2}{1+2a}
		\end{pmatrix}$, $a\neq -\frac{1}{2}$ \\
		\hline
		$E_5(0,y)$	& $\begin{pmatrix}
		0 & 0\\[2mm]
		1 & 0
		\end{pmatrix}, \quad
		\begin{pmatrix}
		0 & 0\\[2mm]
		c_{1,2} & -1
		\end{pmatrix}$, \\
		&
	{} where $c_{1,2}=\frac{-1 \pm \sqrt{1-4y}}{2}$. \\
		\hline
		$E_5(0,y)$ \	& $\begin{pmatrix}
		-1 & 0\\[2mm]
		-1 & -1
		\end{pmatrix}, \quad
		\begin{pmatrix}
		-1 & 0\\[2mm]
		c_{1,2} & 0
		\end{pmatrix}$, \\ [6mm]
	$ y\neq0$	 & 	where $c_{1,2}=\frac{1 \pm \sqrt{1-4y}}{2}$, \ $c_{1,2}\neq 0$ \\
		\hline
		$E_5(0,y)$  &
		$\begin{pmatrix}
		\frac{1-4y+\sqrt{1-4y}}{8y-2} & -\frac{1}{\sqrt{1-4y}}\\[2mm]
		\frac{y}{\sqrt{1-4y}} & \frac{1-4y-\sqrt{1-4y}}{8y-2}
		\end{pmatrix}$, \\ [6mm]
	 $y\neq 0$, $y\neq \frac{1}{4}$	& $\begin{pmatrix}
		\frac{1-4y-\sqrt{1-4y}}{8y-2} & \frac{1}{\sqrt{1-4y}}\\[2mm]
		-\frac{y}{\sqrt{1-4y}} & \frac{1-4y+\sqrt{1-4y}}{8y-2}
		\end{pmatrix}$  \\
        {}&{} \\
		\hline
%	\end{tabular}
%
% \begin{tabular}{|c|l|}
%	\hline
%	Evolution Algebra & Matrices of RBOs of weight 1 \\
%	{} & on  the evolution algebra \\
%		\hline \hline
			{}	& $\begin{pmatrix}
			0 & 1\\[2mm]
			0 & 0
			\end{pmatrix}, \quad
			\begin{pmatrix}
			0 & b_{1,2}\\[2mm]
			0 & -1
			\end{pmatrix}$,\\ [6mm]
			$E_5(x,0)$	& $\begin{pmatrix}
			-1 & -b_{1,2}\\[2mm]
			0 & 0
			\end{pmatrix}, \qquad
			\begin{pmatrix}
			-1 & -1\\[2mm]
			0 & -1
			\end{pmatrix}$,\\
			&
			where $b_{1,2}=\frac{1 \pm \sqrt{1-4x}}{2}$, \ $b_{1,2}\neq 0$. \\
            {}&{} \\
			\hline
			$E_5(x,0)$	&
			$\begin{pmatrix}
			\frac{1-4x+\sqrt{1-4x}}{8x-2} & -\frac{x}{\sqrt{1-4x}}\\[2mm]
			\frac{1}{\sqrt{1-4x}} & \frac{1-4x-\sqrt{1-4x}}{8x-2}
			\end{pmatrix}$, \\ [6mm]
		 $x\neq0$, $x\neq \frac{1}{4}$	& $\begin{pmatrix}
			\frac{1-4x-\sqrt{1-4x}}{8x-2} & \frac{x}{\sqrt{1-4x}}\\[2mm]
			-\frac{1}{\sqrt{1-4x}} & \frac{1-4x+\sqrt{1-4x}}{8x-2}
			\end{pmatrix}$ \\
			\hline
			$E_5(0,0)$	& $\begin{pmatrix}
			-1 & 0\\[2mm]
			0 & 0
			\end{pmatrix}$ \\
			\hline
				$E_5(x,y)$	&  $\begin{pmatrix}
			-1 & 0\\[2mm]
			0 & -1
			\end{pmatrix}$\\
			\hline
			$E_5(x,1-x)$ 	&  $\begin{pmatrix}
			\frac{-3-i \sqrt{3}}{6} & -\frac{i}{\sqrt{3}}\\[2mm]
			\frac{i}{\sqrt{3}} & \frac{-3+i \sqrt{3}}{6}
			\end{pmatrix}$,\\ [6mm]
		$x\neq \frac{1\pm i\sqrt{3}}{2}$	& $\begin{pmatrix}
			\frac{-3+i \sqrt{3}}{6} & \frac{i}{\sqrt{3}}\\[2mm]
			-\frac{i}{\sqrt{3}} & \frac{-3-i \sqrt{3}}{6}
			\end{pmatrix}$ \\ [6mm]
			\hline
				$E_5(x,y)$ & $\begin{pmatrix}
			-1-d & -\frac{d(1+d)}{c}\\[2mm]
			c & d
			\end{pmatrix}$,\\
			$x=\frac{d(1+d)(c+2cd-d(1+d))}{c^2(1+3d+3d^2)}$, &  $d\neq 0, \ d\neq -1, \ d\neq -\frac{3\pm i \sqrt{3}}{6}$, \\
			$y=\frac{c(1-c+2d)}{c^2(1+3d+3d^2)}$	
			& $c\neq 0, \ c\neq \frac{d(1+d)}{1+2d}, \ c\neq 1+2d$.\\
			\hline
%\end{tabular}
%	
% \begin{tabular}{|c|l|}
%	\hline
%	Evolution Algebra & Matrices of RBOs of weight 1 \\
%	{} & on  the evolution algebra \\
%	\hline \hline
	$E_6(0)$	& $\begin{pmatrix}
	\frac{-3+i\sqrt{3}}{6} & 0\\[2mm]
	0 & \frac{-3-i\sqrt{3}}{6}
	\end{pmatrix}, \quad
	\begin{pmatrix}
	\frac{-3-i\sqrt{3}}{6} & 0\\[2mm]
	0 & \frac{-3+i\sqrt{3}}{6}
	\end{pmatrix}$, \\ [7mm]
	&
	$\begin{pmatrix}
	\frac{-3+i\sqrt{3}}{6} & -\frac{i}{\sqrt{3}}\\[2mm]
	\frac{i}{\sqrt{3}} & \frac{-3-i\sqrt{3}}{6}
	\end{pmatrix},\quad
	\begin{pmatrix}
	\frac{-3-i\sqrt{3}}{6} & \frac{i}{\sqrt{3}}\\[2mm]
	-\frac{i}{\sqrt{3}} & \frac{-3+i\sqrt{3}}{6}
	\end{pmatrix}$, \\ [7mm]
	& $\begin{pmatrix}
	\frac{-3-i\sqrt{3}}{6} & -\frac{\sqrt[6]{-1}}{\sqrt{3}}\\[2mm]
	\frac{(\sqrt[6]{-1})^5}{\sqrt{3}} & \frac{-3+i\sqrt{3}}{6}
	\end{pmatrix}$, \quad
	$\begin{pmatrix}
	\frac{-3+i\sqrt{3}}{6} & \frac{\sqrt[6]{-1}}{\sqrt{3}}\\[2mm]
	-\frac{(\sqrt[6]{-1})^5}{\sqrt{3}} & \frac{-3-i\sqrt{3}}{6}
	\end{pmatrix}$,\\ [7mm]
	&$\begin{pmatrix}
	\frac{-3+i\sqrt{3}}{6} & \frac{(\sqrt[6]{-1})^5}{\sqrt{3}}\\[4mm]
	\frac{\sqrt[6]{-1}}{\sqrt{3}} & \frac{-3-i\sqrt{3}}{6}
	\end{pmatrix},\quad
	\begin{pmatrix}
	\frac{-3-i\sqrt{3}}{6} & -\frac{(\sqrt[6]{-1})^5}{\sqrt{3}}\\[4mm]
	\frac{\sqrt[6]{-1}}{\sqrt{3}} & \frac{-3+i\sqrt{3}}{6}
	\end{pmatrix}$.\\ [7mm]
	\hline
	$E_6(x)$	& $\begin{pmatrix}
	-1 & 0\\[2mm]
	0 & -1
	\end{pmatrix}$ \\ [7mm]
	\hline
	$E_6(\frac{-b^3-c^3}{bc^2})$	& $\begin{pmatrix}
	\frac{b^2-c}{2c} & b\\[2mm]
	c & \frac{-b^2-c}{2c}
	\end{pmatrix}$ \\
 $-b^3-c^3\neq 0$	& where the parameters $b,c$ are solutions \\
$b\neq 0, \ c\neq 0$	&   of the following system of equations\\
	& 	$\begin{cases}
	\frac{b^6+5b^3c^3+4c^6}{c}&=\frac{c(b^3+c^3)}{b}, \\[1mm]
	\quad \frac{b^4}{c}+4bc^2&=c.
	\end{cases}$ \\ [7mm]
	\hline
\end{longtable}	
	
\end{thm}

\section*{Acknowledgements}

 We sincerely acknowledge Professor U.A. Rozikov for helpful discussions.
 
This work was partially supported by Agencia Estatal de Investigaci\'on (Spain), grant MTM2016-79661-P and by Xunta 
de Galicia, grant ED431C 2019/10 (European FEDER support included, UE).

\end{document}